%% file: main.tex
\documentclass[11pt, twocolumn]{article}
\input{preamble}

\usepackage{verbatim}
\usepackage{authblk}
\usepackage{wrapfig}
\usepackage{indentfirst}
\usepackage[
    left=0.5in, 
    right=0.5in, 
    top=0.5in, 
    bottom=0.87in
]{geometry}

\usepackage{sectsty}
\sectionfont{\fontsize{11}{15}\selectfont}
\subsectionfont{\fontsize{11}{15}\selectfont}

\usepackage[compact]{titlesec}
\titlespacing{\section}{0pt}{11pt}{5.5pt}
\titlespacing{\subsection}{0pt}{5.5pt}{5.5pt}
\titleformat{\section}
  {\sffamily\bfseries}
  {\thesection.}{1em}{}
\titleformat{\subsection}
  {\sffamily\it}
  {\thesubsection}{1em}{}

\usepackage[
    format=plain, 
    font=it
]{caption}

\makeatletter
\renewcommand*{\@fnsymbol}[1]{\ifcase#1\or\#\else\@arabic{#1}\fi}
\makeatother

\usepackage[
natbib=true,
backend=biber,
style=numeric,
sorting=none,
]{biblatex}
\title{\textbf{
Scalable and Interactive Electricity Grid Expansion Planning\thanks{
This is a paper for the Applied Energy Symposium: MIT A+B, Aug.\ 12-15, 2024, Cambridge, USA.}
}}
\author[1, *]{Anthony Degleris}
\author[1]{Abbas El Gamal}
\author[2]{Ram Rajagopal}

\affil[1]{Department of Electrical Engineering, Stanford University}
\affil[2]{Department of Civil Engineering, Stanford University}
\affil[*]{\textit{(Corresponding Author: degleris@stanford.edu)}}

\date{}

\begin{document}
\maketitle

\renewcommand{\abstractname}{ABSTRACT}
\begin{abstract}
Large scale grid expansion planning studies are essential to rapidly and efficiently decarbonizing the electricity sector.
These studies help policy makers and grid participants understand which renewable generation, storage, and transmission assets should be built and where they will be most cost effective or have the highest emissions impact.
However, these studies are often either too computationally expensive to run repeatedly or too coarsely modeled to give actionable decision information.
In this study, we present an \textit{implicit gradient descent} algorithm to solve expansion planning studies at scale, i.e., problems with many scenarios and large network models.
Our algorithm is also \textit{interactive}: given a base plan, planners can modify assumptions and data then quickly receive an updated plan.
This allows the planner to study expansion outcomes for a wide variety of technology cost, weather, and electrification assumptions.
We demonstrate the scalability of our tool, solving a case with over a hundred million variables.
Then, we show that using warm starts can speed up subsequent runs by as much as 100x.
We highlight how this can be used to quickly conduct storage cost uncertainty analysis.
\end{abstract}


\begin{figure}[t!]
    \centering
    \includegraphics[width=3.65in]{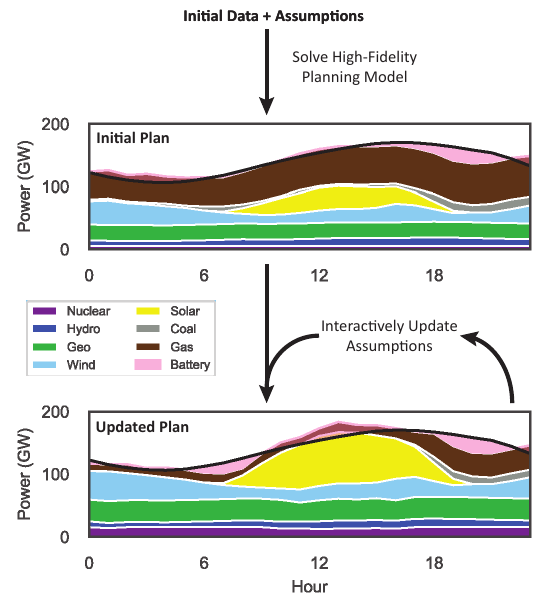}
    \caption{
    Workflow for interactive planning studies.
    The planner first establishes an initial plan on a large scale, high-fidelity grid model.
    Based on initial results, the planner interactively updates assumptions (e.g., an emissions target or the cost of transmission upgrades) and quickly receives a new plan, iteratively repeating this process.
    }
    \label{fig:diagram}
\end{figure}

\section{INTRODUCTION}

Electrification of heating, transportation, and industrial processes is essential to successful decarbonization.
This transformation will increase electric load in the U.S.\ by upwards of 50\%~\citep{National-Renewable-Energy-Laboratory2022-ma} and in turn require significant expansion of electric grid transmission, storage, and generation assets.
Capacity expansion planning studies, which identify high-quality investment plans in electric grid resources, will therefore be key to finding cost-effective paths to decarbonization.

In modern electricity systems, modeling uncertainty is a core a part of an effective planning process.
Specifically, the weather dependence of renewable generation introduces new sources of \textit{short-term} uncertainty in energy supply.
Unpredictable costs associated with new technologies (e.g., grid scale storage) and rapidly shifting policy targets also create a great deal of \textit{long-term} uncertainty about the fundamental assumptions and inputs used in planning models.
In this environment, planners need tools that are both \textit{scalable}, to model many short-term scenarios simultaneously, and \textit{interactive}, to quickly explore and understand the impact of different assumptions on grid outcomes.

In this work, we consider a general bilevel planning problem, called \textit{multi-value expansion planning}.
This problem, in general, is nonconvex, so we describe an efficient, scalable approximate solution method that can be interactively updated post hoc.
Our model, described in Section~\ref{sec:model}, is quite general and allows for jointly planning generation, transmission, storage for a variety of objectives.
Inspired by its remarkable success in solving complex machine learning problems, we use \textit{stochastic gradient descent}, described in Section~\ref{sec:gradient}, to solve multi-value planning problems.
This approach has a simple interpretation: at each iteration, the planner assesses the current expansion plan and uses sensitivity analysis to slightly adjust improve upon the plan.
This method is scalable to high-fidelity grid models with tens of millions of variables and constraints and can be interactively updated, i.e., the planner can change input data or assumptions and quickly find a new plan.
In Section~\ref{sec:results}, we test the performance of our method on a variety of Western U.S. test cases and demonstrate applications to sensitivity analysis and exploratory planning.
We conclude in Section~\ref{sec:conclusion} and suggest several practical use cases for our method.

\subsection{Related Work}

Existing expansion planning tools fall primarily into one of two categories: (1) simplified models of the grid~\citep{Fishbone1981-yx, Loulou2008-wz, Fripp2012-ih} that can be studied for many different scenarios and assumptions, and (2) high-fidelity models~\citep{Majidi-Qadikolai2018-hb, Neumann2019-mj} that accurately model grid physics (e.g., DC linearized power flow) but are too computationally intensive to run repeatedly for differing assumptions and data.
Most models also assume a central planner that jointly optimizes operations and investment to minimize total cost, indirectly incorporating emissions through carbon taxes and caps;
accurately modeling decentralized electricity markets requires a \textit{bilevel expansion planning} problem \citep{Pozo2017-kz, Gonzalez-Romero2020-zp} that significantly increases computational complexity.
Finally, planning tools are \textit{static} by default, i.e., they produce a single plan that cannot be modified.
Several studies try to address this limitation using sensitivity analysis~\citep{Yue2018-oa, Sanajaoba-Singh2018-ed, Shirizadeh2022-jn} and by exploring the space of near-optimal solutions, i.e., modeling to generate alternatives~\citep{DeCarolis2016-bl, Price2017-mg, Neumann2021-hp}.
For a comprehensive review of expansion planning tools, we refer the reader to~\citet{Gonzalez-Romero2020-zp}.

\section{MATERIALS AND METHODS}

\subsection{Model}
\label{sec:model}

We consider an electricity network parametrized by a vector $\eta \in \reals^K$, which includes information about the constructed generator, transmission, and storage capacity throughout the network.
The \textit{multi-value expansion planning problem} is to choose the network capacities to minimize investment costs plus some function of the operational outcomes,
\begin{equation}
\begin{array}{ll} \label{eq:mvp}
    \textrm{minimize} \quad
    & \gamma^T \eta + \frac{1}{S} \sum_{s=1}^S h_s(x_s^*(\eta))  \\[0.5em]
    \textrm{subject to}
    & \eta^{\min} \leq \eta \leq \eta^{\max},  
\end{array}
\end{equation}
where the variable is the network parameters $\eta \in \reals^K$.
The functions $h_s : \reals^N \rightarrow \reals$ are called the \textit{planner's objectives} and may include things such as the total cost of operation, total grid emissions, or some function of the locational marginal prices.
The function $x_s^*(\eta) : \reals^K \rightarrow \reals^N$ is called the \textit{dispatch map} for scenario $s$,
\begin{equation} \label{eq:dispatch}
    x_s^*(\eta) = \underset{x}{\argmin} \ \Big( c_s(x)\ \textrm{ s.t. }\ A_s(\eta) x \leq b_s(\eta) \Big), 
\end{equation}
which models short-term operational decisions under scenario $s \in \{1, \ldots, S\}$.
The dispatch map is simply the outcome of the economic dispatch problem or production cost model used by the planner.
The functions $c_s : \reals^N \rightarrow \reals$ model the cost of operations (e.g., fuel cost) for that scenario.
Each scenario has $M$ constraints $A(\eta) x \leq b(\eta)$, where $A(\eta) \in \reals^{M \times N}$ and $b(\eta) \in \reals^M$ are functions of the electricity network capacities $\eta$.
For example, if $\eta_1$ is the capacity of a generator, then one of the constraints will include that generator's power limits, which depend on $\eta_1$.
Individual scenarios can represent an hour, a day, or even a week of operations, and each scenario models either (a) a different period in time or (b) a different set of assumptions about future technology costs, load growth, etc.

Problem~\eqref{eq:mvp} is a \textit{bilevel optimization problem}, i.e., an optimization problem that depends on the outcome of another problem.
When $h_s = c_s$, multi-value planning can be recast as a joint optimization problem over $\eta$ and $x_1, \ldots, x_S$;
this corresponds to classical expansion planning, which minimizes total long-term system cost for a vertically integrated utility.
When $h_s \neq c_s$, we can consider more complex problems where the planner's objectives do not necessarily align with that of a cost-based market clearing mechanism.
In this case, the problem becomes significantly more complex and requires more sophisticated algorithms to solve.
We give an example of a multi-value problem below.

\textbf{Emissions-aware planning.}
In many situations, the planner (a grid operator or state policy maker) may want to decrease emissions through target  investment in new transmission, renewable generation, and storage.
However, in the absence of a sufficient carbon tax, decentralized markets will always clear on monetary cost.
The planner can still reduce emissions and account for the market-clearing equilibrium by setting $h$ to include a weighted emissions penalty.


\subsection{Implicit Gradient Descent}
\label{sec:gradient}

Traditionally, solutions to bilevel expansion planning problems are found by rewriting~\eqref{eq:mvp} as a single level problem using the optimality conditions of the inner problem.
The resulting problem is a nonconvex quadratic program and is difficult to solve exactly in a reasonable time.
Standard approaches include applying off-the-shelf non-linear solvers, branch-and-bound algorithms, and disjunctive reformulations (for integer problem); see~\citet{Pozo2017-kz, Gonzalez-Romero2020-zp, Wogrin2020-xp} for further details.

We propose forgoing the single level reformulation propose and instead applying \textit{gradient descent} directly to~\eqref{eq:mvp}.
At each iteration $i$, we compute the gradient of the objective,
\begin{equation*}
    \Delta^{\iter}
    \coloneq 
    \gamma + \frac{1}{S} \sum_s (\partial x^*_s(\eta^{\iter}))^T \cdot \nabla h_s(x^*(\eta^{\iter})).
\end{equation*}
For large problems, we use stochastic estimates of the gradient by sampling $B$ scenarios and averaging them to further reduce compute time;
in this case, the algorithm is called \textit{stochastic gradient descent}.
Then, the planner updates the current network parameters $\eta^\iter$ to
\begin{equation*}
    \eta^{\iterplus} 
    \coloneq 
    \underset{ [\eta^{\min}, \eta^{\max}] }{\mathbf{proj}} \left( \eta^\iter - \alpha \Delta^{\iter} \right),
\end{equation*}
where the parameter $\alpha > 0$ is the step size, $\nabla h_s : \reals^N \rightarrow \reals^N$ is the gradient of $h_s$, and $\partial x^*_s : \reals^K \rightarrow \reals^{N \times K}$ is the Jacobian of $x^*_s(\eta)$.
Computing the Jacobian of $x^*_s(\eta)$ requires differentiating an $\argmin$ operator.
We achieve this using the \textit{implicit function theorem}~\citep{Dontchev2014-ck} on the optimality conditions of~\eqref{eq:dispatch}.
This technique, called \textit{implicit differentiation}, has been successfully used to apply machine learning and optimization to various physical systems; refer to~\citet{Blondel2021-iu} and references therein for further details.
Further details for the specific application of this approach to multi-value planning problems can be found in~\citet{Degleris2024-oq}.

\textbf{Interpretability.}
Gradient descent has a natural interpretation as sensitivity-based planning.
Given the currently planned network capacities $\eta^{\iter}$, the planner computes the sensitivity of different operational outcomes with respect to the parameters, $\partial x_s^*(\eta^{\iter})$.
The planner uses these sensitivity to slightly improve the planned network capacities, then repeats this process iteratively.

\textbf{Scalability.}
Although gradient descent is only guaranteed to converge to a local stationary point, we empirically observe near-optimal performance in many cases~\citep{Degleris2024-oq}.
Moreover, gradient descent scales well to large problems ($K$, $N$, and $M$) with many scenarios ($S$).
Specifically, each iteration of the algorithm requires solving the dispatch model for every scenario $s \in \{1, \ldots, S\}$ and computing its dispatch map Jacobian $\partial x_s^*(\eta)$;
this can be parallelized across all $S$ scenarios.
However, the number of iterations to converge does not depend on the size of the problem or the number of scenarios~\citep{Degleris2024-oq}.
In practice, the algorithm converges in a few dozen to a few hundred iterations.

\textbf{Interactivity.}
The initial network capacities $\eta^{(0)}$ can be chosen a variety of ways:
they can be set to the existing network capacities ($\eta^{(0)} = \eta^{\min}$), set heuristically, or set using the outcome of a previous planning study.
As we will see in Section~\ref{sec:experiment-warm-start}, the third approach, called \textit{warm starting}, significantly reduces the solve time and enables \textit{interactive} planning studies, in which the planner decides on an initial plan, then incrementally tweaks the plan based on past results.

\section{RESULTS}
\label{sec:results}

In this section, we first demonstrate the need for high-fidelity grid models by showing that naive scenario selection (e.g., peak day planning) underperforms when evaluating on a full year of data.
We then show that warm starts can be used to accelerate gradient descent by as much as 100x, converging in just a few dozen iterations.
Finally, we use warm starts to rapidly explore planning outcomes for different battery storage costs.

We implement implicit gradient descent in Python for generic multi-value planning problems using \verb|cvxpy|~\citep{Diamond2016-ja} and Mosek~\citep{MOSEK-ApS2022-cf} to solve the dispatch model $x_s^*(\eta)$ and PyTorch~\citep{Paszke2019-yo} to differentiate $h_s$.
We develop a custom, modular modeling system for computing the Jacobian $\partial x_s^*(\eta)$ via the implicit function theorem that is simple to extend to different grid devices.

All our experiments use the PyPSA-USA Western Interconnect dataset~\citep{Tehranchi2023-jo}.
Specifically, we configure cases of varying sizes (100 to 500 nodes) and use 2019 weather data coupled with 2050 forecasts for load and technology costs to produce 8760 hours of data from NREL~\citep{National-Renewable-Energy-Laboratory2022-ma}.
We solve the emissions-aware planning problems described in Section~\ref{sec:model} with varying carbon weights, jointly planning generation, storage, AC line, and DC line capacities.
We treat each day (with 24 hourly time periods) as an independent scenario, allowing us to study problems with between 1 and 365 scenarios.
We run all experiments on a single AMD EPYC 7763 processor with 64 cores (128 threads) at NERSC.

\begin{figure}
    \centering
    \includegraphics{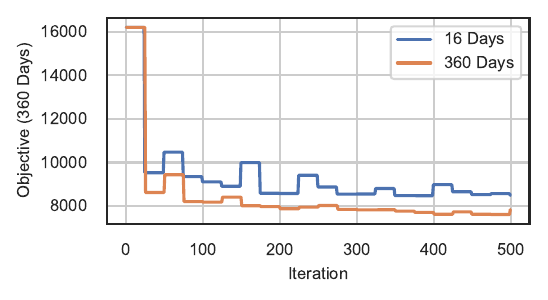}
    \caption{
    Objective values (evaluated on 360 days of data) for stochastic gradient descent solved on 16 days of data (blue curve) and 360 days of data (orange curve).
    Both runs use a batch size $B=8$, take the same time to complete, and are evaluated on the full 360 days.
    Sampling from the full dataset leads to 11.3\% reduction in objective value.
    }
    \label{fig:scalability}
\end{figure}

\subsection{Scalability}

We first study the scalability of our method and demonstrate that high-fidelity models find better quality solutions than their approximate counterparts.
Specifically, we solve two multi-value problems with a carbon weight of \$200 / ton CO$_2$: one with 16 key days\footnote{
We select the top 4 peak load, peak net load, peak renewable, and lowest renewable days as our 16 key days.}
and another with 360 scenarios (roughly a full year).
We run stochastic gradient descent with a batch size $B=8$ on both problems for 500 iterations and evaluate their performance on the full 360 scenarios every 25 iterations.

Each scenario considers 24 hours of operation for a network with 500 nodes, 1799 generators, 1191 AC lines, 3 HVDC lines, and 576 batteries.
We plan the capacity of all 3569 devices simultaneously across all scenarios of hourly data.
If formulated as a joint problem (as is done traditionally), the full optimization problem would have over a hundred million variables and a similar number of constraints.

We plot the 360 scenario objective value at each iteration for both cases in Fig.~\ref{fig:scalability}.
Solving the 360 scenario case directly leads to a 11.6\% lower objective value, suggesting that higher resolution models can help identify more efficient investments. 
Critically, both algorithms have the same batch size and therefore the same runtime, approximately 85 seconds per iteration.
About half of this time is spent evaluating the objective on the full dataset, which can be performed less frequently to save time.

\begin{figure}[h!]
    \centering
    \includegraphics{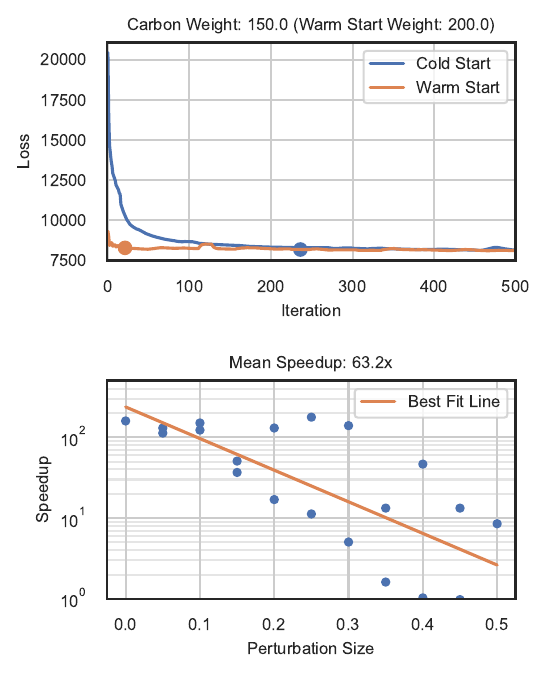}
    \caption{
    Warm starts significantly speedup solve times for large planning problems.
    (Top) Convergence plots for 500 iterations of stochastic gradient descent applied to a emissions-aware planning problem with a \$150.0 / ton CO$_2$ carbon weight.
    Blue curve: loss when initialized from "no-expansion", which converges in 235 iterations.
    Orange curve: loss when warm started from the solution to a similar problem with carbon weight \$200.0 / ton CO$_2$, which converges in just 21 iterations.
    (Bottom) Speedup achieved by warm starting as a function of the size of the perturbation to the carbon weight.
    For small perturbations, warm starting leads to upwards of a 100x speedup.
    }
    \label{fig:warm-start-performance}
\end{figure}

\begin{figure*}[t]
\centering
\includegraphics{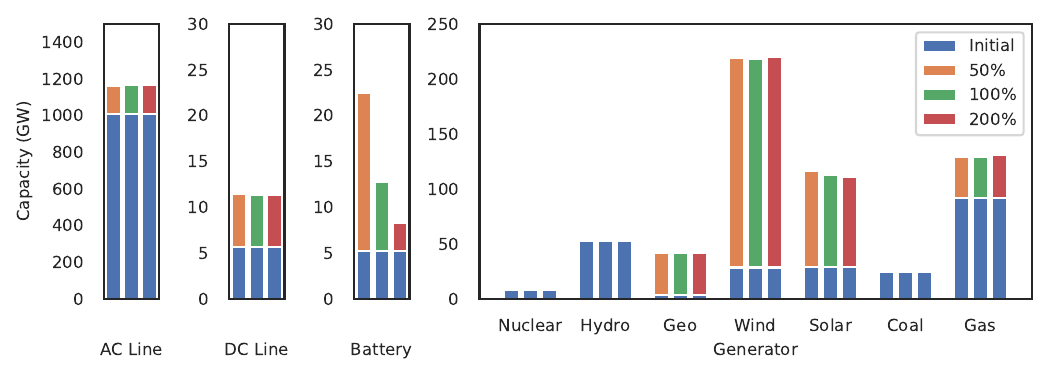}
\caption{
Warm starting enables rapid exploration of potential expansions under different battery cost assumptions.
(Blue bars) Initial capacity of each asset before expansion.
(Orange, green, and red bars) Expanded capacity of each asset battery capital costs are 50\%, 100\%, and 200\%, respectively, of NREL ATB~\citep{National-Renewable-Energy-Laboratory2022-ma} forecasts.
Optimal investment in transmission and generation capacity are stable across different battery cost assumptions;
optimal investment in battery storage depends on assumed technology costs.
}
\label{fig:battery-expansion}
\end{figure*}

\subsection{Warm Start Performance}
\label{sec:experiment-warm-start}

In this experiment, we analyze how warm starting can be used to reduce compute time when perturbing problem inputs.
We consider 16 scenario, 100 node cases for different carbon weights ranging from \$100 to \$300 / ton CO$_2$.
For each carbon weight, we solve the problem with two different initializations: (cold start) no expansion, i.e., all capacities are set to their existing values, and (warm start) the solution to the same problem with a carbon weight of \$200 / ton CO$_2$.
Each problem is solved using using 500 iterations stochastic gradient descent with a batch size $B=4$; a single iteration takes 6 to 10 seconds.
To compare the cold start and warm start solutions, we say each algorithm converges when its loss is within 2\% of the best cold start loss.

Results are displayed in Fig.~\ref{fig:warm-start-performance}.
We first show the loss plots for the carbon weight of \$150~/~ton CO$_2$.
Despite the modest (25\%) perturbation to the carbon weight, the warm start run still converges more than 10x faster than the cold start run.
Next, we measure the speedup (warm start iterations divided by cold start iterations) for different perturbation sizes.
On average, warm starts lead to a 63.3x speedup, with upwards of a 100x speedup for small perturbation sizes.

In practical terms, this means a planning study that might normally take about 30 minutes (250 iterations) when run from scratch might only take 30 seconds (5 iterations) to update to for new parameters and input data.
This could allow planners to interactively explore how proposed plans changed under different assumptions and inputs.

\subsection{Analyzing Battery Cost Assumptions}

Using the same configuration as Section~\ref{sec:experiment-warm-start}, with a carbon weight of \$200 / ton CO$_2$, we use warm starts to rapidly explore planning results for 20 different values assumed battery investment costs between 10\% and 200\% of the NREL ATB forecast.
We display the resulting expansions for three different assumptions on battery investment cost in Fig.~\ref{fig:battery-expansion}.
As expected intuitively, battery expansion decreases significantly as a function of battery investment cost.
However, we find that investment in all other grid assets essentially does not vary for different battery costs.
During an interactive study, a planner could use information like this to determine minimal regret investments in transmission in generation; regardless of the ultimate cost of battery storage, a similar amount of transmission and renewable generation will be needed.

\section{CONCLUSION}
\label{sec:conclusion}

In this work, we describe a scalable, interactive algorithm for complex expansion planning problems.
Our algorithm, which applies gradient descent to the planning problem using implicit differentiation, can solve high-fidelity grid expansion planning problems in reasonable amounts of time.
The algorithm can also leverage warm starts to rapidly modify an existing plan to new data and assumptions, and has a simple interpretation as iterative sensitivity analysis.

We believe our tool could be of great interest to grid planners and policy makers who could use such a tool to interactively explore future expansion plans for different assumptions, goals, and data on high-fidelity grid models.
For example, policy makers could explore the implications of different tax levels or emissions targets without forgoing complex constraints like generator ramping limits or linearized power flow, and planners could investigate the impact of climate and electrification assumptions on capacity requirements.
Tools like this could be critical in designing coherent, robust, and cost effective plans for upgrading electricity systems.

\section*{ACKNOWLEDGEMENTS}

This research used resources of the National Energy Research Scientific Computing Center (NERSC),
a U.S. Department of Energy Office of Science User Facility located at Lawrence Berkeley National Laboratory, operated under Contract No. DE-AC02-05CH11231 using NERSC award DDRERCAP0026889.

This material is based upon work supported by the U.S. Department of Energy, Office of Science, Office of Advanced Scientific Computing Research, Department of Energy Computational Science Graduate Fellowship under Award Number DE-SC0021110

\printbibliography[title=REFERENCES]


\end{document}

%% file: preamble.tex
\usepackage{mathtools} 
\usepackage{xcolor} 
\usepackage{tikz} 
\usetikzlibrary{shapes.geometric, positioning, arrows}

\newcommand{\reals}{\mathbf{R}}

\newcommand{\argmin}{\mathrm{argmin}}

\newcommand{\iter}{{(i)}}
\newcommand{\iterplus}{{(i+1)}}